\documentclass[12pt,oneside,reqno]{amsart}
\usepackage[margin=1in]{geometry}
\usepackage{color}
\usepackage{tikz-cd}
\usepackage{esint,amssymb}
\usepackage{graphicx}
\usepackage{MnSymbol}
\usepackage{mathtools}
\usepackage[colorlinks=true, pdfstartview=FitV, linkcolor=blue, citecolor=blue, urlcolor=blue,pagebackref=false]{hyperref}
\usepackage{microtype}
\usepackage{amsmath}
\usepackage{xifthen}
\usepackage{verbatim}
\definecolor{darkgreen}{rgb}{0,0.5,0}
\definecolor{darkblue}{rgb}{0,0,0.7}
\definecolor{darkred}{rgb}{0.9,0.1,0.1}
\usepackage[all]{xy}
\usepackage{enumitem} 
\usepackage{tikz}
\usepackage{xcolor}
\usepackage{hyperref}
\definecolor{myfuchsia}{RGB}{250,0,220}
\definecolor{mypurple}{RGB}{140,0,200}
\hypersetup{colorlinks=true, citecolor=myfuchsia, linkcolor=mypurple}

\newtheorem{theorem}{Theorem}
\newtheorem{proposition}[theorem]{Proposition}
\newtheorem{lemma}[theorem]{Lemma}
\newtheorem{corollary}[theorem]{Corollary}

\theoremstyle{definition}
\newtheorem{remark}[theorem]{Remark}

\newtheorem{definition}[theorem]{Definition}
\newtheorem{conjecture}[theorem]{Conjecture}

\newcommand{\cref}[1]{Corollary~\ref{c.#1}}

\numberwithin{equation}{section}
\numberwithin{theorem}{section}

\newcommand{\N}{\mathbb{N}}
\newcommand{\R}{\mathbb{R}}

\newcommand{\C}{\mathbb{C}}

\newcommand{\T}{\mathsf{T}}
\renewcommand{\subset}{\subseteq}

\newcommand{\test}[1][]{%
\ifthenelse{\equal{#1}{}}{omitted}{given}%
}

\renewcommand{\bar}{\overline}

\renewcommand{\part}{\partial}

\usepackage{dsfont}
\usepackage{slashed}

\newcommand{\SL}{\operatorname{SL}}
\newcommand{\Tr}{\operatorname{Tr}}

\newcommand{\Id}{\operatorname{Id}}

\begin{document}
\title{A counterexample on the main flatness conjecture for $\SL_n(\C)$-Fock bundles}

\author[Sotiria Chatzimarkou]{Sotiria Chatzimarkou}
\address[Sotiria Chatzimarkou]{University of Patras}
\email{sotiria.mathematics@gmail.com}
\subjclass[2010]{}
\date{\today}
\begin{abstract}
The main flatness conjecture for $\SL_n(\C)$-Fock bundles predicts the existence of a unique strongly positive flat compatible Hermitian metric for every $\SL_n(\C)$-Fock bundle. We use cyclic Higgs bundles to disprove this conjecture for every rank $n \geq 3.$
\end{abstract}

\maketitle
\tableofcontents
\section{Introduction}
Let $\Sigma$ be a closed, oriented Riemannian surface of genus $g$ at least $2$. In his foundational paper \cite{hitchin1987self}, Hitchin established a remarkable connection between Teichm\"uller theory and the theory of Higgs bundles. More precisely, he proved that the Teichm\"uller space associated to $\Sigma$, namely the moduli space that parametrizes complex structures on $\Sigma$ modulo diffeomorphisms isotopic to the identiry, is diffeomorphic to the Hitchin component of the character variety $\mathcal{X}(\pi_1(\Sigma),\mathrm{PSL}_2(\R))$ via the Non-Abelian Hodge correspondence. Motivated by this identification, Hitchin in \cite{hitchin1992lie} introduced the first family of higher Teichm\"uller spaces, the Hitchin components using the theory of Higgs bundles. Labourie \cite{labourie2004anosov} and independently Fock and Goncharov \cite{fock2004moduli} proved that they are indeed higher Teichm\"uller spaces. 

Much later, Fock and Thomas in \cite{fock2021higher, Thomas2022} introduced higher complex structures generalizing the classical notion of complex structures and conjectured the existence of a mapping class group equivariant diffeomorphism between the moduli space of higher complex structures and the $\mathrm{PSL}_n(\R)$-Hitchin component, analogous to the Non-Abelian Hodge correspondence in rank two. The conjecture was proved in rank three by Nolte in \cite{nolte2022canonical}, using the positive answer to Labourie's conjecture for rank three stated in \cite{labourie2008cross}. However, its validity in higher rank remains a challenging problem.

A different approach to understanding higher complex structures is provided by gauge theory. In \cite{kydonakis2023}, Kydonakis, Reid and Thomas introduced Fock bundles, which provide a gauge-theoretic realization of higher complex structures. Based on this realization they stated the following conjecture:
\begin{conjecture}(The main flatness conjecture) \cite{kydonakis2023} [Conjecture 4.14]
    Let $(E,b,v,\Phi_0)$ an $\SL_n(\C)$-Fock bundle with a Hermitian metric $h_0$ compatible with $b$ and let $\rho_0=\rho_{h_0}$. Then there \textbf{exists} a \textbf{unique} Hermitian endomorphism-valued function $\eta \in \Omega^0(\Sigma,\mathfrak{sl}(E)^{\sigma_b,-\rho_0})$ such that $\Phi=e^{-\eta}\Phi_0 e^\eta$ is positive and the corresponding 3-term connection $\Phi+D_{\Phi,h_0}+\Phi^{*_{h_0}}$ is flat, meaning that it satisfies the equation 
\begin{equation}
    \mathcal{M}_\Phi (h_0)=0.
\end{equation}  
\end{conjecture}

If true, this conjecture would imply the Fock-Thomas conjecture. Kydonakis, Reid and Thomas only established the main flatness conjecture locally in a neighborhood of the Fuchsian locus. The validity of the conjecture globally remained an open problem in \cite{kydonakis2023}. 

In this paper, we show that the main flatness conjecture is false. More precisely, we construct a counterexample in every rank $n\geq3$ showing that there exists an explicit Fock bundle where either existence or uniqueness fails.

Our approach uses cyclic Higgs bundles. Let $(E,\theta_t,Q_t)$ be a cyclic Higgs bundle with $\theta_t=\Phi_F+\mathcal{Q}_t, \hspace{0.4cm} \mathcal{Q}_t=tq_nE_{1n}, \hspace{0.2cm}t\geq 0,$ and $h_t$ its harmonic metric. Consider the Fock field $\Phi_t=\Phi_F+\mathcal{Q}_t^{*_{h_t}}$. The associated flat connection becomes the 3-term canonical connection that is used in \cite{kydonakis2023} as long as $h_t$ is positive for $\Phi_t$. 

There is a global Fock regrouping such that the two components of $\Phi_t=\Phi_F+\mathcal{Q}_t^{*_{h_t}}$ are $A_t=\sum\limits_{i=1}^{n-1} a_{i,t}E_{i+1,i}$ and $B_t=\beta_tE_{n1}.$ The norm of the graph operator $[A_t,X]\rightarrow [B_t,X]$ is given by $\|C_{\Phi_t,h_t}\|_{op}^2=|\beta_t|^2\sum \limits _{i=1}^{n-1}|a_{i,t}|^{-2}.$ We prove that $h_t$ is strongly positive if and only if $\|C_{\Phi,h,x}\|^2_{op}<1.$ 

The asymptotic estimates for harmonic metrics of cyclic Higgs bundles obtained by Collier and Li \cite{collier2014asymptotics} imply that, away from the zeros of $q_n$ we have $|\frac{\beta_t}{a_{i,t}}|\rightarrow 1$ as $t\to \infty$. Therefore, $\mathcal{W}_t(x):=\|C_{\Phi_t,h_t}\|_{op}^2\rightarrow n-1.$ Let $\Lambda(t)=\max \limits_{x\in \Sigma} \mathcal{W}_t(x)$, then for $n\geq 3$ and away from the zeros of $q_n$, $\Lambda(t)>1$ for $t$ large. Then, for $t_*=\inf\{t>0:\Lambda(t)\geq 1\}$, $\Lambda(t_*)=1$ since $\Lambda$ is continuous. We call $t_*$ the positivity wall. We use the existence of this positivity wall to prove the following
\begin{theorem}
For every rank $n\geq3$ the main flatness conjecture \cite[Conjecture 4.14]{kydonakis2023} is false. Specifically, at least one of the following conditions is forced at the positivity wall:
\begin{enumerate}[label=\roman*)]
\item the field $\Phi_{t_*}$ has no strongly positive flat compatible Hermitian metric $h_{t_*}$,
\item for all $t<t_*$ sufficiently close to $t_*$, the field $\Phi_t$ has at least two distinct strongly positive flat compatible Hermitian metrics $h_t, \hat{h}_t$.
\end{enumerate}
\end{theorem}
This theorem shows that universal existence and universal uniqueness cannot both hold in rank $n\geq3$. Indeed, if no strongly positive flat compatible Hermitian metric exists at the wall, the existence part of the conjecture fails. On the other hand, if a strongly positive flat compatible Hermitian metric exists, then there is another strongly positive flat compatible Hermitian metric exactly before the wall and hence uniqueness fails. Therefore, the assertion that every $\SL_n(\C)$-Fock bundle admits a unique strongly positive flat compatible Hermitian metric cannot hold in general. This establishes the failure of the main flatness conjecture in every rank $n\geq 3.$

\subsection*{Acknowledgments} I would like to thank Panagiotis Dimakis for various fruitful conversations throughout the course of this project and in particular for his help with the analytic arguments in the paper. I would like to thank the Brin Mathematics Research Center for giving me the opportunity to spend a month at the University of Maryland, during which much of this paper was written. 

\subsection*{AI declaration} The free version of ChatGPT(OpenAI) was used solely for proofreading purposes, including spelling, grammar and syntactic corrections. 

\section{Preliminaries}
\subsection{Basic geometric notions}
Let $\Sigma$ be a closed, connected, oriented, smooth surface of genus $g \geq 2$ equipped with a smooth complex vector bundle $E$. 

\begin{definition}
    A complex structure on $\Sigma$ is an automorphism $J:T^*\Sigma\longrightarrow T^*\Sigma$ such that $J^2=-\Id$.
\end{definition} 
Given a complex structure J, consider the Riemann surface X=($\Sigma$,J). Let $K=T^{*1,0}X$ denote the canonical bundle and $\bar{K}=T^{*0,1}X$ its dual.

An $\SL_n (\C)$-structure on $E$ is a fixed determinant trivialization $v:\det E \rightarrow \C$, $v$ is a complex volume form on $\Sigma$.

Let $\mathfrak{sl}(E)$ be the set of trace-free endomorphisms of $E$ equipped with the invariant trace pairing $B_{\Tr} (U,V)=\Tr(UV)$. 

Consider a complex-bilinear symmetric form $b:E\oplus E\longrightarrow \C $ such that $\det_v (b)=1$. The transpose with respect to $b$ is given by
\begin{equation} \label{b}
b(Xu,v)=b(u,X^{\dagger_b}v). 
\end{equation}
We define the involution 
\begin{equation}
    \sigma_b(X)=-X^{\dagger _b}.
\end{equation}
The equation $\sigma_b(X)=-X$ is equivalent to $X=X^{\dagger_b}$ or, in words, that $X$ is $b$-self adjoint. Given a Hermitian metric $h$ the adjoint $X^{\dagger_h}$ defines another involution 
\begin{equation}
    \rho_h(X)=-X^{\dagger_h}. 
\end{equation}
For $a \in \Omega ^p(\mathfrak{sl}(E))$ and $\beta\in \Omega ^q(\mathfrak{sl}(E))$ their commutator is given by: 
\begin{equation}
    [a,\beta](u,v)=[a(u),\beta(v)]-[a(v),\beta(u)] 
\end{equation}
and 
\begin{equation}
[a,a](u,v)=2[a(u),a(v)].
\end{equation}

Next we define cyclic Higgs bundles and state their properties required for disproving the main flatness conjecture. 

\begin{definition}
An $\SL_n(\C)$-cyclic Higgs bundle over $\Sigma$ is a Higgs bundle $(E,\Phi)$ admitting a holomorphic splitting 
\[
E\cong L_1\oplus...\oplus L_n,~ \bigotimes\limits_{i=1}^n L_i \cong \mathcal O_X
\]
such that 
\begin{equation*}
        \Phi =\begin{pmatrix}
            0 & 0 & ... & 0 &q_n \\
            q_1 &0 & ... & 0 &0 \\
            0 &q_2 & ... &0 &0 \\
            \vdots & \ddots & \ddots & \vdots &\vdots \\
            0 &0 & ... &q_{n-1} & 0,
        \end{pmatrix} 
\end{equation*}
where $q_i \in H^0(\Sigma, K\otimes L_i^{-1}\otimes L_{i+1})$ and $L_{n+1} = L_1$.    
\end{definition}

For $\SL(n,\C)$-cyclic Higgs bundles the harmonic metric splits with respect to the holomorphic splitting 
\[
h = \mathrm{diag}(h_1,...,h_n),~ \prod\limits_{i=1}^n h_i = 1.
\]
For a cyclic Higgs bundle in the Hitchin section, $q_2=...=q_{n-1} =1$ and thus the Higgs bundle is specified by $q_n\in \Gamma(K^n)$. The metric satisfies $h_{n+1-i} = h_i^{-1}$.

\subsection{Fock bundles and their properties}
In this subsection we introduce the notion of Fock bundles and discuss their properties. We start by defining principal nilpotent elements. 
\begin{definition}
An element $N \in\mathbf{sl}_n(\C) $ is principal nilpotent if it is nilpotent and has a single Jordan block. This is equivalent to $N^n=0$ and $N^{n-1}\neq 0$, which is equivalent to its centralizer having dimension $n-1$.
\end{definition}
The standard principal nilpotent element is:
\begin{equation}
     S_n=\sum \limits _{j=1}^{n-1} E_{j+1,j}.
\end{equation}
\begin{remark}\label{remark}
     If $A \in \mathfrak{sl}_n(\C)$ is a principal nilpotent element then its centralizer $Z_{\mathfrak{sl}_n} (A)=\mathrm{span}_\C\{A,A^2,...,A^{n-1}\}.$ Evidently, $c_1A+...+c_{n-1}A^{n-1}$ is principal nilpotent if and only if $c_1\neq0$.
\end{remark}
\begin{definition}\label{fock bundle}
    Given $(E,v)$, an $\SL_n(\C)$-Fock bundle is a quadruple $(E,b,v,\Phi)$ where $\Phi \in \Omega^1(\Sigma,\mathfrak{sl}(E))$ and $b$ is a complex bilinear symmetric form, defined as in \eqref{b}, such that:
\begin{enumerate}[label=\roman*)]
  \item $[\Phi \wedge \Phi]=0,$
  \item $\Phi_x (v)$ is principal nilpotent for every $x \in \Sigma$ and for every non-zero $v \in T_x \Sigma,$ 
  \item $\sigma_b(\Phi)=-\Phi$.
\end{enumerate}
The 1-form $\Phi$ is called a Fock field.
\end{definition}
\begin{remark}
    Let $u,v \in T_x\Sigma$ then $\Phi_x(u), \Phi_x(v)$ are principal nilpotent elements and $[\Phi_x(u),\Phi_x(v)]=0$ by the axiom (iii) of Definition \ref{fock bundle}.
\end{remark}
We give a proof for the following Proposition stated in \cite{kydonakis2023}[Proposition 3.2]
\begin{proposition}\label{structure}
Every $\SL_n(\C)$-Fock bundle induces a complex structure on $\Sigma$.
\end{proposition}
 \begin{proof}
     Let $(E,v,b,\Phi)$ an $\SL_n(\C)$-Fock bundle. For every $x \in \Sigma$ we prove that the set 
     \begin{equation}\label{complex line}
         L_{\Phi,x}:=\{\xi \in T_c^\C\Sigma:\Phi_c^\C(\xi) \text{ is not principal nilpotent}\}
         \end{equation}
         is a complex line that contains no nonzero real vector, the lines vary smoothly with $x$ and 
         \begin{equation} \label{athrisma}
         T^\C\Sigma=L_\Phi \oplus \bar{L_\Phi}.
         \end{equation}
         Then by defining 
         \begin{equation}
         \begin{split}
         T_{J_\Phi}^{0,1}\Sigma:=L_\Phi, \\T_{J_\Phi}^{1,0}\Sigma:=\bar{L_\Phi}
         \end{split}
         \end{equation}
         we obtain a complex structure on $\Sigma$.

     Fix $x \in \Sigma$ and a real basis $u,v$ of $T_x\Sigma$. Set $A=\Phi_x(u)$ and $B=\Phi_x(v)$ that are principal nilpotent elements and $[A,B]=0$ by the previous remark. By remark \ref{remark} there are unique coefficients $c_1,...,c_{n-1} \in \C$ such that 
    \begin{equation}\label{B}
         B=c_1A+c_2A^2+...+c_{n-1}A^{n-1}.
    \end{equation}
B is principal nilpotent thus $c_1\neq 0.$
Let $\xi=zu+wv \in T_x^\C\Sigma$ then
\begin{equation}
    \Phi_x^\C(\xi)=(z+wc_1)A+wc_2A^2+...+wc_{n-1}A^{n-1}. 
\end{equation}
By remark \ref{remark} this element is principal nilpotent if and only if $z+wc_1\neq 0$. Therefore, 
\begin{equation}\label{line} 
L_{\Phi,x}=\C(v-c_1u)
\end{equation}
which proves the existence of the complex line \ref{complex line}.
If $L_{\Phi,x}$ contains a nonzero real vector then there are $0\neq z,w \in \R$ such that $z+wc_1=0$. Then $\Phi_x(zu+wv)$ is not principal nilpotent which contradicts axiom (ii) of Definition \ref{fock bundle}. Thus 
\begin{equation}\label{tomi}
    L_{\Phi,x}\cap T_x\Sigma={0}.
\end{equation}
For a complex line in $T^\C\Sigma$, equation \eqref{tomi} is equivalent to $L_{\Phi,x}\neq \bar{L_{\Phi,x}}$, which implies \eqref{athrisma}.
Choose $u,v$ smooth real vector fields forming a local frame then $A=\Phi(u)$ is a smooth family of principal nilpotent elements and $A,A^2,...,A^{n-1}$ are linearly independent. Thus, they form a smooth local frame of the centralizer $Z(A)\subset \mathfrak{sl}(E)$. If we express $B=\Phi(v)$ in this local frame we obtain that every $c_j$ in \eqref{B} is smooth. Thus, by \eqref{line}, $L_\Phi$ has a smooth local generator. 
There is a smooth real bundle endomorphism $J_\Phi$ such that $J_\Phi^2=-I$ that acts by $-i$ on $L_\Phi$ and by $i$ on $\bar{L_\Phi}$. 

 \end{proof}

 By Proposition \ref{structure}, given a local holomorphic coordinate $z$ we can locally write $\Phi=\Phi_1 dz+\Phi_2d\bar{z}$. The first condition of definition \ref{fock bundle} implies that $[\Phi_1, \Phi_2]=0$. 

 \begin{lemma}
 Fix a complex coordinate on $\Sigma$. An $\SL_n(\C)$-Fock bundle is locally of the form $\Phi=\Phi_1 dz+\Phi_2d\bar{z}$ where:
 \begin{equation*}
     \Phi_1=\sum \limits _{j=1} ^{n-1}E_{j+1,j}= 
     \begin{pmatrix} 
     0 & 0 & ... & 0 & 0 \\ 
     1 & 0 & ... & 0 & 0 \\ 
     0 & 1 & ... & 0 & 0 \\
    \vdots & \vdots &\ddots & \vdots &\vdots \\ 
    0 & 0 & ... &1 &0
     \end{pmatrix}
 \end{equation*}
 and 
\begin{equation*}
     \Phi_2=\mu_2 \Phi_1 +\mu _3 \Phi _1^2 + ...+\mu_n \Phi _1^{n-1}= \begin{pmatrix} 
     0 & 0 & ... & 0 & 0 \\ 
     \mu_2 & 0 & ... & 0 & 0 \\ 
     \mu_3 & \mu_2 & ... & 0 & 0 \\
    \vdots & \vdots &\ddots & \vdots &\vdots \\ 
    \mu_n & \mu_{n-1} & ... &\mu_2 &0
     \end{pmatrix}
\end{equation*}
where $\mu_k$ are local holomorphic functions such that $\mu_2 \bar{\mu}_2 \neq 1$.
\end{lemma}
\begin{proof}
     $\Phi_1$ is principal nilpotent so we choose a local holomorphic frame in which $\Phi_1=S_n$.
     $[\Phi_1,\Phi_2]=0$ implies that $\Phi_2 \in Z(\Phi_1)=\mathrm{span}\{\Phi_1,\Phi_1^2,...,\Phi_1^{n-1}\}$. Therefore, we can uniquely write $\Phi_2=\mu_2 \Phi_1 +\mu _3 \Phi _1^2 + ...+\mu_n \Phi _1^{n-1},$ where $\mu_k$ are local holomorphic functions.
\end{proof}
\begin{definition}
    We call a Hermitian metric $h$ with determinant equal to $1$ compatible with the bilinear symmetric form $b$ \eqref{b} if
    \begin{equation}
        \rho_h \sigma _b=\sigma_b \rho_h.
    \end{equation}
\end{definition}
For a Fock field $\Phi$ we define its $h$-adjoint as
\begin{equation}
    \Phi^{*_h}:=-\rho _h (\Phi).
\end{equation}
If locally $\Phi=\Phi_1 dz+\Phi_2d\bar{z}$ then $\Phi^{*_h}=(\Phi_1)^{\dagger _h} dz+(\Phi_2)^{\dagger_h}d\bar{z}$.
\begin{proposition}
   For every $\SL_n(\C)$-Fock bundle $(E,b,v,\Phi)$ equipped with a hermitian metric $h$ compatible with $b$, $(E,b,v,\Phi^{*_h})$ is also an $\SL_n(\C)$-Fock bundle.
\end{proposition}
\begin{proof}
    One can easily check that the three conditions of the definition of Fock bundles are true
\begin{enumerate}[label=\roman*)]
  \item $[\Phi^{*_h} \wedge \Phi^{*_h}]=0,$
  \item $\Phi^{*_h}_x (v)$ is principal nilpotent for every $x \in \Sigma$ and for every non-zero $v \in T_x \Sigma$ since the set of principal nilpotent matrices is invariant under hermitian conjugation
  \item $\sigma_b(\Phi^{*_h})=-\sigma_b \rho_h(\Phi)=-\rho_h\sigma_b(\Phi)=-\rho_h(-\Phi)=-\Phi^{*_h}$.
\end{enumerate}
\end{proof}

For every point $x \in \Sigma$, consider the vector space $V_x=T_x^{\star} \Sigma\oplus_\R \mathfrak{sl}(E_x)$. We denote
\begin{equation}
    I_{\Phi,x}=\mathrm{Im}([\Phi\wedge*]:\mathfrak{sl}(E_x)\rightarrow V_x)
\end{equation}
and
\begin{equation}
    K_{\Phi,x}=\ker([\Phi\wedge*]:V_x\rightarrow \Lambda^2T_x^*\Sigma\oplus_\R \mathfrak{sl}(E_x)).
\end{equation}
Then, the first condition of the definition of Fock bundles \ref{fock bundle} implies that the following sequence is a complex
\begin{equation}
\mathfrak{sl}(E_x)\xrightarrow{[\Phi\wedge*]}V_x\xrightarrow{[\Phi\wedge*]}
\Lambda^2T_x^*\Sigma\oplus_\R \mathfrak{sl}(E_x).
\end{equation}
\begin{definition}
    Two Fock fields $\Phi, \Psi\in \Omega^1(\Sigma,\mathfrak{sl}(E))$ are called transverse if for every $x \in \Sigma$ 
\begin{equation}
        V_x=K_{\Phi,x}\oplus I_{\Psi,x}=K_{\Psi,x}\oplus I_{\Phi,x}.
\end{equation}
\end{definition}
\vspace{0,5cm}
For a Hermitian metric $h$ compatible with $b$ there exists a pseudo-Hermitian form $q_{h,x}$ on $V_x$. For $\alpha, \beta\in V_x$ 
\[
q(\alpha,\beta):=q_{h,x}(\alpha,\beta) = -\frac{i}{2}\Tr(\alpha^{\star}\wedge\beta). 
\]
Writing $\alpha = A\,dz + B\,d\bar z$ locally, the form is given by 
\begin{equation}\label{pseudo}
    q_{h,x}(A dz+B d\bar{z},A dz+B d\bar{z})=\|A\|^2_{h,x}-\|B\|^2_{h,x}=(\Tr(AA^{\dagger_h})-\Tr(BB^{\dagger_h}))\,dz\,d\bar z.
\end{equation}
\begin{definition}
    The hermitian metric $h$ is strongly positive for $\Phi$ if $q_{h,x}$ is positive definite on $I_{\Phi,x}$ for every $x \in \Sigma$.
\end{definition}
Write the Fock field as $\Phi=\Phi_1 dz+\Phi_2 d\bar{z}$. Since  $\Phi_1$ is principal nilpotent and $\Phi_2$ commutes with $\Phi_1$,  $\ker([\Phi_1 \wedge*])\subseteq \ker([\Phi_2 \wedge*])$. This implies that there is a well-defined map
\begin{equation}
\begin{split}
    C_{\Phi,h,x}:\mathrm{Im}([\Phi_1 \wedge*])&\longrightarrow \mathfrak{sl}(E_x) \\
    [\Phi_1,X] &\longmapsto [\Phi_2,X]
\end{split}
\end{equation}
\begin{lemma} \label{lemma}
The hermitian metric $h$ is strongly positive at $x$ if and only if $\|C_{\Phi,h,x}\|_{\mathrm{op}}^2<1$, where $\|L\|_{\mathrm{op}}=\inf \{c\ge 0: \|Lv\|\le c\|v\|$, for every $v \in V \}$ is the operator norm for a linear map $L:V\rightarrow W$ .
\end{lemma}
\begin{proof}
    Let $\phi=[\Phi_1,X]dz+[\Phi _2,X]d\bar{z}$ an element of $I_{\Phi,x}$. We denote $Y=[\Phi_1,X]$. Then the pseudo-Hermitian norm of $\phi$ is equal to $\|Y\|^2-\|C_{\Phi,h,x}Y\|^2$ which is positive for every $Y\neq 0$ exactly when $\|C_{\Phi,h,x}\|_{op}^2<1$.
\end{proof}
Thus a hermitian metric $h$ is strongly positive exactly when
\begin{equation}
    \max \limits _{x \in \Sigma} \|C_{\Phi,h,x}\|_{op}^2<1.
\end{equation}
Strong positivity implies that $I_\Phi$ and $K_{\Phi^{*_h}}$ are orthogonal complements with respect to the pseudo-Hermitian form \ref{pseudo} which is equivalent to the transversality of $\Phi, \Phi^{*_h}$. 
\begin{theorem}(Filling Theorem)\label{filling thm} Let $(E,b,v,\Phi)$ an $\SL_n(\C)$-Fock bundle equipped with a hermitian metric h compatible with $b$ and $\Phi, \Phi^{*_h}$ are transverse. Then there exists a unique connection $D_{\Phi,h}$ satisfying the equations:
\begin{enumerate}[label=\roman*)]
  \item $D_{\Phi,h}v=0$
  \item $D_{\Phi,h}b=0$
  \item $D_{\Phi,h}h=0$
  \item $D_{\Phi,h}\Phi=D_{\Phi,h}\Phi^{*_h}=0$.
\end{enumerate}
\end{theorem}
We refer the reader to \cite[Theorem 4.6]{kydonakis2023} for the  proof.
\vspace{0.5cm}

According to this Theorem, in every $\SL_n(\C)$-Fock bundle there exists a 3-term family of connections, with middle term $D_{\Phi,h}$. The 3-term family of connections is given by
\begin{equation}\label{fam}
    \nabla(\lambda)=\lambda^{-1}\Phi+D_{\Phi,h}+\lambda\Phi^{*_h},
\end{equation}
where $\lambda \in \C^*$.

Let 
\begin{equation}\label{m}
\mathcal{M}_\Phi (h):=F_{D_{\Phi,h}}+[\Phi \wedge\Phi^{*_h}].
\end{equation}
The curvature of $\nabla(\lambda)$ is equal to $F_{\nabla(\lambda)}=\mathcal{M}_\Phi (h).$
The flatness of the family $\nabla(\lambda)$ is equivalent to
\begin{equation}\label{flat eq}
   \mathcal{M}_\Phi (h)=0
\end{equation}
which is a nonlinear PDE similar to the famous Hitchin equation introduced in \cite{hitchin1987self}.
\begin{definition}
     A Hermitian metric $h$ compatible with $b$ is called strongly positive flat compatible metric if it is strongly positive and solves the equation $\mathcal{M}_\Phi (h)=0$.
\end{definition}

\subsection{Linearization and local invertibility} Before we state the main flatness conjecture, we introduce some further notation and some of the analytic machinery introduced in \cite{kydonakis2023} to study equation \eqref{flat eq} and used in the disproof of the main flatness conjecture below.

Fix now any smooth strongly positive pair $(\Phi,h)$ and write $\Phi^{\star} = \Phi^{\star_h}$, $\rho = \rho_h$ and
\begin{equation}\label{cartan}
\begin{split}
\mathfrak p &= \mathfrak p_{\sigma,\rho} := \{u\in \mathfrak{sl}(E) : \sigma_b(u) = u,~ \rho(u)=-u\},\\
\mathfrak t &= \mathfrak t_{\sigma,\rho} := \{u\in \mathfrak{sl}(E) : \sigma_b(u) = u,~ \rho(u)= u\}. 
\end{split}
\end{equation}

For $\eta \in \Gamma(\mathfrak p)$ define 
\[a_\eta=e^\eta,~ h_\eta(u,v)=h(a_\eta^{-1}u,a_\eta^{-1}v),~\Phi_\eta=e^{-\eta}\Phi e^\eta.\]

Let $\mathfrak w = \mathfrak {sl}(E)^{\sigma_b}$. Since $h$ is assumed to be strongly positive, the pseudo-Hermitian form $q_{h,x}$ defined above is positive definite on $I_{\Phi}$. By compactness of $\Sigma$ there exists $\varepsilon>0$ such that 
\[
q(\alpha,\alpha) \geq \varepsilon |\alpha|_{g,h}^2
\]
for all $\alpha\in I_{\Phi}$. From now on, fix a metric $g$ on $\Sigma$. Let us define 
\[
Q_h(\alpha,\beta) = \int_{\Sigma} q(\alpha,\beta)\,d\mathrm{vol}_g,~ \|\alpha\|^2_{+,h} = Q_h(\alpha,\alpha).
\]
Let $\mathcal I:= I_{\Phi}^{\sigma_b} = I_{\Phi}\cap (T^{\star}\Sigma\otimes_{\R}\mathfrak w)$ and $\mathcal J:= I_{\Phi^{\star}}^{\sigma_b}$. It is proven in \cite{kydonakis2023} that 
\[
T^{\star}\Sigma\otimes_{\R}\mathfrak w = \mathcal I\oplus \mathcal J.
\]
Define $Q_{\Phi,h}:= \Id -2\pi_{\Phi}$,
where $\pi_{\Phi}: T^{\star}\Sigma\otimes_{\R}\mathfrak w \to \mathcal I$ is the projection.  

Let $D = D_{\Phi,h}$ denote the filling connection. Given $\eta\in \Gamma(\mathfrak p)$ consider 
\[
\Phi_s = e^{-s\eta}\Phi e^{s\eta},~ S_s= D_{\Phi_s,h},~ M_s=F_{D_s} + [\Phi_s\wedge\Phi_s^{\star}].
\]
The linearization 
\[
\frac{\,dM_s}{\,ds} = L_{\Phi_s,h} \eta 
\]
maps $L_{\Phi,h}:\Gamma(\mathfrak p) \to \Omega^2(\mathfrak t)$ and is given by the expression
\[
L_{\Phi,h} \eta = d_DQ_{\Phi,h}d_D\eta+[[\Phi,\eta]\wedge\Phi^{*_h}]-[\Phi\wedge[\Phi^{*_h},\eta]].
\]
It is easier for us to work with the operator 
\[
\mathcal A_{\Phi,h}:= -\frac{i}{2}\star_gL_{\Phi,h} :\Gamma(\mathfrak p) \to \Gamma(\mathfrak p),
\]
where $\star_g$ is the Hodge star operator for the metric $g$. It satisfies the identity 
\begin{equation}
\begin{split}
\langle \eta, \mathcal A_{\Phi,h}\zeta\rangle_{L^2} &= -\frac{i}{2}\int_{\Sigma} \Tr (\eta L_{\Phi,h}\zeta)\\
&= 2\mathrm{Re}( Q_h(\pi_{\Phi} \,d_D\eta,\pi_{\phi}\,d_D\zeta)) + 2\mathrm{Re}( Q_h([\Phi,\eta],[\Phi,\zeta])).
\end{split}
\end{equation}
The right hand side is symmetric in $\eta,\zeta$ so $\mathcal A$ is formally self-adjoint. Furthermore, since $\,d_D\eta$ is $\mathfrak p$-valued 
\[
2q(\pi_{\Phi} \,d_D\eta,\pi_{\phi}\,d_D\eta) \geq \frac{\varepsilon}{2}|\,d_D\eta|^2_{g,h},
\]
and since $u\to [\Phi,u]$ is injective on each fiber of $\mathfrak p$, 
\[
q([\Phi,\eta],[\Phi,\eta]) \geq \mu |u|_h^2
\]
for $u\in \mathfrak p$ by continuity and compactness. These together imply that 
\begin{equation}\label{zerokernel}
c\|\eta\|^2_{H^1} \leq -\frac{i}{2}\int _\Sigma \Tr(\eta L_{\Phi,h}\eta) ~\text{for}~c>0.
\end{equation}
Finally, the symbol $\sigma(\mathcal A)$ of $\mathcal A$ satisfies the identity
\[
\langle u,\sigma(\mathcal A)(x,\xi)u\rangle_u = 2q(\pi_{\Phi}(\xi u),\pi_{\Phi}(\xi u)) \geq \frac{\varepsilon}{2}|\xi|_g^2|u|_h^2. 
\]
Therefore, $\mathcal A_{\Phi,h}: H^{k+2}(\mathfrak p)\to H^k(\mathfrak p)$ for $k\geq 0$ is a uniformly elliptic, formally self-adjoint operator. Inequality \eqref{zerokernel} implies that the kernel is trivial and formal self-adjoint-ness implies that as a map between these Sobolev space $\mathcal A$ is an isomorphism. 
 
Finally, we are in a position to state and prove the proposition required for deriving the contradiction in the main theorem of the paper. 

\begin{proposition}\label{local existence}
Let $t\mapsto\Phi_t$ a smooth family of Fock fields on a fixed datum $(E,\nu,b)$. If $h_0$ is a strongly positive flat Hermitian metric at $t_0$ then for $t$ sufficiently close to $t_0$ there exists a unique smooth family $h(t)$ of compatible strongly positive flat hermitian metrics with $h(t_0) = h_0$. 
\end{proposition}
\begin{proof}
Choose $k\geq 2$. On a neighborhood of $(t_0,0)\in \R\times H^{k+2}(\mathfrak p_{\sigma,\rho_0})$ let
\[
\mathcal F(t,\eta) = \mathcal M_{e^{-\eta}\Phi_te^{\eta}}(h_0)\in H^k(\Omega^2(\mathfrak t_{\sigma,\rho_0}))
\]
The choice of $k$ guarantees that this is a smooth map between the given Banach spaces. Also, it holds that 
\[
D_{\eta} \mathcal F(t_0,0) = L_{\Phi_{t_0},h_0}
\]
which is an isomorphism between the given Banach spaces. Thus the Banach implicit function theorem gives a unique smooth local map $t\to \eta(t)$ with $\eta(t_0) = 0$ and $\mathcal F(t,\eta(t))=0$. Since $k\geq 2$, the Sobolev embedding theorem implies that $H^{k+2}(\mathfrak p)$ embeds continuously into $C^0$. Since strong positivity is an open condition we can choose a small enough neighborhood around $t_0$ so that all $\eta(t)$ give rise to strongly positive Hermitian metrics. This finishes the proof.
\end{proof}

\subsection{The main flatness conjecture for $\SL_n(\C)$-Fock bundles}
The following conjecture, stated in \cite[Conjecture 4.14]{kydonakis2023} is similar to the Non-Abelian Hodge correspondence and claims that for every $\SL_n(\C)$-Fock bundle there exists a unique strongly positive flat compatible Hermitian metric $h$ making the family of connections \eqref{fam} flat. The purpose of this paper is to prove that this conjecture fails to be true for every rank $n\ge 3$.

\begin{conjecture}(The main flatness conjecture) Let $(E,b,v,\Phi_0)$ an $\SL_n(\C)$-Fock bundle with a Hermitian metric $h_0$ compatible with $b$ and let $\rho_0=\rho_{h_0}$. Then there exists a unique Hermitian endomorphism-valued function $\eta \in \Omega^0(\Sigma,\mathfrak{p}_{\sigma_b,\rho_0})$ such that $\Phi=e^{-\eta}\Phi_0 e^\eta$ is positive and the corresponding 3-term connection $\Phi+D_{\Phi,h_0}+\Phi^{*_{h_0}}$ is flat, meaning that it satisfies the equation 
\begin{equation}
    \mathcal{M}_\Phi (h_0)=0.
\end{equation}  
\end{conjecture}

\subsection{From cyclic Higgs bundles to Fock bundles}

Finally, we describe the procedure for obtaining a Fock field from a Higgs field. 

For a closed Riemann surface $X=(\Sigma,J)$ and $K$ the canonical bundle, fix a theta characteristic $\kappa$ such that $\kappa^2\cong K$ and let $(E,\theta)$ be an $\SL_n(\C)$-cyclic Higgs bundle. Let $b_F$ be the canonical anti-diagonal holomorphic symmetric form $b_F$ associated with this cyclic Higgs bundle given by pairing $L_i$ with $L_{n+1-i}$ such that $det(b_F)=1.$ For $q_n \in H^0(X,K^n)$ and $t\ge 0$, write
\begin{equation}
    \mathcal{Q}_t=tq_nE_{1n}
\end{equation}
and consider the Higgs field
\begin{equation}
    \theta_t=\Phi_F+\mathcal{Q}_t.
\end{equation}
The Higgs bundle $(E,\theta_t)$ is stable for every $t\ge 0$ and has a unique harmonic metric $h_t$ with determinant $1$ of the form 
\begin{equation}
    h_t=\mathrm{diag}({h_{1,t},h_{2,t},...,h_{n,t})}
\end{equation}
such that $h_{j,t}h_{n+1-j, t}=1$. This metric is also compatible with $b_F$. Let $D_t$ be the Chern connection of the metric $h_t$. Define the Fock field 
\begin{equation}\label{fock field}
    \Phi_t=\Phi_F+\mathcal{Q}_t^{*_{h_t}}.
\end{equation}
\begin{proposition} \label{2.21}
For every $n\ge 3$ and for every $t\ge0$, $(E_X,v,b_F,\Phi_t)$ is an $\SL_n(\C)$-Fock bundle inducing the complex structure $J$. Also 
\begin{equation}
    \Phi_t^{*_{h_t}}=\Phi_F^{*_{h_t}}+\mathcal{Q}_t
\end{equation}
and
\begin{equation}
    \Phi_t+\Phi_t^{*_{h_t}}=\theta_t+\theta_t^{*_{h_t}}.
\end{equation}
\end{proposition}
\begin{proof}
   Let 
   \begin{equation*}
L_j=\kappa^{n+1-2j}
\end{equation*}
and
\begin{equation*}
    E_X=\bigoplus_{j=1}^{n}L_j
\end{equation*} 
and fix a holomorphic coordinate. We want to show that $\Phi_t$ is a Fock field. 
\begin{enumerate}[label=\roman*)]
    \item Locally we have
\begin{equation}\label{ab}
\begin{split}
    &\Phi_t=Adz+Bd\bar z, \\
    &A_t=\sum\limits_{j=1}^{n-1}a_{j,t}E_{j+1,j}, \hspace{0,3cm} a_{j,t}>0 \\
    &B_t=\beta_tE_{n1}
    \end{split}
\end{equation}
Since $A_t^{n-1}=(a_{1,t}...a_{n-1,t})E_{n1}$, we have $B_t=\mu_tA_t^{n-1}$, for $\mu_t=\frac{\beta_t}{a_{1,t}...a_{n-1,t}}$. Thus, we can easily see that $[A_t,B_t]=0$ and $[\Phi_t\wedge\Phi_t]=0$.
\item Let $v\in T_x\Sigma$ a non zero real vector and $\xi=dz(v)$ then $\xi\neq 0$ and 
\begin{equation*}
    \Phi_t(v)=\xi A_t+\bar{\xi}\mu_tA_t^{n-1}=\xi A_t\big(I+\frac{\bar{\xi}}{\xi}\mu_tA_t^{n-2}\big).
\end{equation*}
The element $A_t^{n-2}$ is nilpotent thus the factor $I+\frac{\bar{\xi}}{\xi}\mu_tA_t^{n-2}$ is invertible and commutes with $A_t$.Therefore, $\Phi_t(v)$ is principal nilpotent by \ref{remark}.
\item We have $\sigma_{b_F}(\Phi_F)=-\Phi_F$ and $\sigma_{b_F}(\mathcal{Q}_t)=-\mathcal{Q}_t$. Compatibility implies that $\sigma_{b_F}(\Phi_t)=-\Phi_t$.
\end{enumerate}
This proves that $\Phi_t$ is a Fock field making $(E_X,v,b_F,\Phi_t)$ an $\SL_n(\C)$-Fock bundle.
\end{proof}

\begin{proposition}\label{connection}
    The Chern connection $D_t$ preserves $v, b_F$ and $h_t$ and satisfies the equations
    \begin{equation}
    \begin{split}\label{e}
        & d_{D_t}\Phi_t=0 \\
        & d_{D_t}\Phi_t^{*_{h_t}}=0.
    \end{split}
    \end{equation}
Moreover, the family of connections 
\begin{equation}\label{ee}
    \nabla_t=D_t+\Phi_t+\Phi_t^{*_{h_t}}
=D_t+\theta_t+\theta_t^{*_{h_t}}
\end{equation}
is flat.
\end{proposition}
\begin{proof}
   Since $b_F$ is a holomorphic form, the Chern connection preserves it. The Hermitian metric $h_t$ is compatible with $b_F$ thus the Chern connection also preserves $h_t$. 
   Both $\Phi_F$ and $\mathcal{Q}_t$ are holomorphic K-valued endomorphisms, thus their (0,1)-part vanish by holomorphicity and their (1,0)-part vanish on a complex structure, which implies that
   \begin{equation*}
       d_{D_t}\Phi_F=d_{D_t}\mathcal{Q}_t=0.
   \end{equation*}
   The same is true for their adjoints which proves \eqref{e}.
$D_t+\Phi_t+\Phi_t^{*_{h_t}}$ is flat since it is the harmonic metric equation. If $h_t$ is strongly positive then $\Phi_t$ and $\Phi_{t}^{*_{h_t}}$ are transverse. Thus, by Theorem \ref{filling thm} the uniqueness condition identifies $D_t$ with the connection that makes the family flat.
\end{proof}

When the cyclic harmonic metric $h_t$ is strongly positive we call $D_t$ the canonical connection and the associated 3-term connection $D_t+\Phi_t+\Phi_t^{*_{h_t}}$ is called the canonical 3-term connection.

\begin{remark}
    At a positivity wall, that we will construct in the next section, equations \eqref{e} and \eqref{ee} still hold and give a smooth flat 3-term decomposition. However, the connection $D_t$ is not the canonical 3-term connection that is used in \cite{kydonakis2023} because transversality may fail and the uniqueness condition of Theorem \ref{filling thm} does not hold.
\end{remark}

\section{Disproof of the main flatness Conjecture }

\subsection{A weighted commutator norm}
In what follows we equip $\mathfrak{sl}_n(\C)$ with the usual Hilbert-Schmidt norm. The proof of our main theorem crucially relies on the following
\begin{theorem}\label{wln}
    Let 
    \begin{equation*}
        A=\sum_{i=1}^{n-1}a_iE_{i+1,i}
    \end{equation*}
    for $a_i\neq 0$
    and 
    \begin{equation*}
        B=\beta E_{n,1}.
    \end{equation*}
     Then the map given by
     \begin{equation}
     \begin{split}
         C_{A,B}:\mathrm{Im}([A\wedge*])&\longrightarrow \mathfrak{sl}_n(\C)\\
         [A,X]&\longmapsto[B,X]
    \end{split}
     \end{equation}
    is  well defined and satisfies the identity
     \begin{equation}\label{id}
         \|C_{A,B}\|^2_{op}=|\beta|^2\sum _{i=1}^{n-1}|a_i|^{-2}
     \end{equation}
     and it is attained by a trace-free diagonal matrix X.
\end{theorem}
\begin{proof}
    One can easily see that 
    \begin{equation*}
        B=\frac{\beta}{a_1...a_{n-1}}A^{n-1}.
    \end{equation*}
    Thus, every matrix X that commutes with A also commutes with B and that implies that $C_{A,B}$ is well defined.
    For the proof of the identity let $Y=[A,X]$ and $Z=[B,X]$. 
    \begin{equation*}
       A= \begin{pmatrix}
              0 & 0 & ... &0 & 0\\
        a_1 &0 &...&0 &0 \\
        0 & a_2 &...&0 &0 \\
        \vdots & \vdots & \ddots &\vdots & \vdots \\
        0 &0 &... &a_{n-1} &0
        \end{pmatrix}
    \end{equation*}
    \begin{equation*}
       B= \begin{pmatrix}
             0 & 0 & ... &0 & 0\\
        0 &0 &...&0 &0 \\
        0 & 0 &...&0 &0 \\
        \vdots & \vdots & \ddots &\vdots & \vdots \\
        \beta &0 &... &0 &0
        \end{pmatrix}
    \end{equation*}
    We calculate:
    \begin{equation} \label{1}
        Z_{n,j}=-\frac{\beta}{a_{j-1}}Y_{1,j-1}, \hspace{1cm} 2\leq j\leq n
    \end{equation}
    \begin{equation}\label{2}
        Z_{i,1}=-\frac{\beta}{a_i}Y_{i+1,n}, \hspace{1cm} 1\leq i\leq n-1
    \end{equation}
    \begin{equation}\label{3}
        Z_{n,1}=\beta \sum_{i=1}^{n-1}\frac{Y_{i+1,i}}{a_i}.
    \end{equation}
    Let $M=|\beta|^2 \sum\limits _{i=1}^{n-1}|a_i|^{-2}$ then each coefficient $-\frac{\beta}{a_i}$ of the equations \ref{1} and  \ref{2} has square $|-\frac{\beta}{a_i}|^2\leq M$. Thus, by squaring \ref{1} and \ref{2} we obtain 
    \begin{equation*}
    \begin{split}
    |Z_{n,j}|^2=\left|-\frac{\beta}{a_{j-1}}\right|^2|Y_{1,j-1}|^2\leq M|Y_{1,j-1}|^2 \\
        |Z_{i,1}|^2=\left|-\frac{\beta}{a_i}\right|^2|Y_{i+1,n}|^2\leq M|Y_{i+1,n}|^2
        \end{split}
    \end{equation*}
    By the Cauchy-Schwarz inequality applied to the equation \ref{3} we obtain
    \begin{equation*}
    |Z_{n,1}|^2=|b|^2 \left|\sum\limits_{i=1}^{n-1} \frac{Y_{i+1,i}}{a_i}\right|^2\leq |b|^2 \left(\sum\limits_{i=1}^{n-1}|{a_i}|^{-2}\right)\left(\sum\limits_{i=1}^{n-1} |Y_{i+1,i}|^2\right)\leq M\sum\limits_{i=1}^{n-1} |Y_{i+1,i}\big|^2.
    \end{equation*}
    By summing these three equations we obtain $\|Z\|^2\leq M\|Y\|^2.$ Which implies the identity \eqref{id}.
To attain equality let $X=\mathrm{diag}(x_1,...,x_n)$ such that $x_i-x_{i+1}=c|a_i|^{-2}$ for $1\leq i\leq n-1$. 
 We calculate:
    \begin{equation*}
\begin{split}
   Y&= \left[\begin{pmatrix}
        0 & 0 & ... &0 & 0\\
        a_1 &0 &...&0 &0 \\
        0 & a_2 &...&0 &0 \\
        \vdots & \vdots & \ddots &\vdots & \vdots \\
        0 &0 &... &a_{n-1} &0
    \end{pmatrix} ,
    \begin{pmatrix}
        x_1 & 0 & ... &0 & 0 \\
        0 &x_2 &... &0 &0 \\
        \vdots & \vdots & \ddots &\vdots & \vdots \\
        0 & 0 & ... &x_{n-1} &0 \\
        0 & 0 &... & 0 &x_n
    \end{pmatrix}\right]= \\
   &= \begin{pmatrix} 
    0 & 0 & ... &0 & 0 \\
    a_1x_1 &0 & ... &0 & 0 \\
    0 &a_2x_2 &... &0 &0  \\
    \vdots & \vdots & \ddots &\vdots & \vdots \\
   0 &0 &...& a_{n-1}x_{n-1} &0
    \end{pmatrix} -
    \begin{pmatrix}
         0 & 0 & ... &0 & 0 \\
         a_1x_2 &0 &... & 0 &0 \\
         0 &a_2x_3 &... &0 &0  \\
         \vdots & \vdots & \ddots &\vdots & \vdots \\
        0 & 0 & ... &x_{n-1} &0 \\
        0 &0 &...& a_{n-1}x_{n} &0
    \end{pmatrix}= \\
    &=\begin{pmatrix}
     0 & 0 & ... &0 & 0 \\
    a_1(x_1-x_2) &0 & ... &0 & 0 \\
    0 &a_2(x_2-x_3) &... &0 &0  \\
    \vdots & \vdots & \ddots &\vdots & \vdots \\
   0 &0 &...& a_{n-1}(x_{n-1}-x_n) &0   
    \end{pmatrix}
    \end{split}
\end{equation*}
    and 
    \begin{equation*}
    \begin{split}
        Z&=\left[\begin{pmatrix}
            0 & 0 & ... &0 & 0\\
        0 &0 &...&0 &0 \\
        0 & 0 &...&0 &0 \\
        \vdots & \vdots & \ddots &\vdots & \vdots \\
        \beta &0 &... &0 &0
        \end{pmatrix},
        \begin{pmatrix}
        x_1 & 0 & ... &0 & 0 \\
        0 &x_2 &... &0 &0 \\
        \vdots & \vdots & \ddots &\vdots & \vdots \\
        0 & 0 & ... &x_{n-1} &0 \\
        0 & 0 &... & 0 &x_n
        \end{pmatrix}\right] \\
        &=\begin{pmatrix}
            0 & 0 & ... &0 & 0\\
        0 &0 &...&0 &0 \\
        0 & 0 &...&0 &0 \\
        \vdots & \vdots & \ddots &\vdots & \vdots \\
        \beta x_1 &0 &... &0 &0
        \end{pmatrix}- \begin{pmatrix}
            0 & 0 & ... &0 & 0\\
        0 &0 &...&0 &0 \\
        0 & 0 &...&0 &0 \\
        \vdots & \vdots & \ddots &\vdots & \vdots \\
        \beta x_n &0 &... &0 &0
        \end{pmatrix} \\
        &= \begin{pmatrix}
            0 & 0 & ... &0 & 0\\
        0 &0 &...&0 &0 \\
        0 & 0 &...&0 &0 \\
        \vdots & \vdots & \ddots &\vdots & \vdots \\
        \beta(x_1-x_n) &0 &... &0 &0
        \end{pmatrix}
        \end{split}
    \end{equation*}
One can easily see that $Y_{1,j-1}=0$, $Y_{i+1,n}=0$ and $Y_{i+1,i}\neq0$. So the equations \ref{1} and \ref{2} are equal to 0 and so it remains to show equality only for \ref{3}. By Cauchy-Schwartz the equality holds and that completes the proof.
\end{proof}

\begin{corollary}
    Let 
    \begin{equation*}
        S_n=\sum\limits_{i=1}^{n-1}E_{i+1,i}
    \end{equation*}
    such that $S_n^{n-1}=E_{n,1}$. Then 
    \begin{equation}\label{sup}
        \sup\limits_{[S_n,X]\neq0}\frac{\|[S_n^{n-1},X]\|^2}{\|[S_n,X]\|^2}=n-1.
    \end{equation}
\end{corollary}
\begin{proof}
    In Theorem \ref{wln} set every $a_i=1$ and $\beta=1$ and we obtain
    \begin{equation*}
        A=\begin{pmatrix}
           0 & 0 & ... &0 & 0\\
        1&0 &...&0 &0 \\
        0 & 1 &...&0 &0 \\
        \vdots & \vdots & \ddots &\vdots & \vdots \\
        0 &0 &... &1 &0  
        \end{pmatrix}=S_n
    \end{equation*}
    and 
    \begin{equation*}
        B=\begin{pmatrix}
             0 & 0 & ... &0 & 0\\
        0 &0 &...&0 &0 \\
        0 & 0 &...&0 &0 \\
        \vdots & \vdots & \ddots &\vdots & \vdots \\
        1 &0 &... &0 &0
        \end{pmatrix}=S_n^{n-1}.
    \end{equation*}
Thus, $M=n-1$ and by the Cauchy-Schwarz inequality $\|Z\|^2\leq M\|Y\|^2$ we obtain that $\frac{\|[S_n^{n-1},X]\|^2}{\|[S_n,X]\|^2}=\frac{\|Z\|^2}{\|Y\|^2}\leq n-1$ which proves \ref{sup}. 
\end{proof}
\subsection{The posivity wall}
\subsubsection{The wall function}
In what follows let $\Phi_t=\Phi_F+\mathcal{Q}_t^{*_{h_t}}$ the Fock field we defined in \ref{fock field}. For a local coordinate let $q_n=f(dz)^n$ and let $h_t=\mathrm{diag}(h_{1,t},...,h_{n,t})$ such that $h_{j,t}h_{n+1-j,t}=1$. In Proposition \ref{2.21} the coefficients of \eqref{ab} satisfy the equations
\begin{equation}\label{*}
    |a_{i,t}|=\frac{h_{i+1,t}}{h_{i,t}} , \hspace{0.4cm} |\beta_t|^2=t^2|f|^2\frac{h_{1,t}}{h_{n,t}}.
\end{equation}
We define the wall function 
\begin{equation}\label{w}
    \mathcal{W}_t(x):=\|C_{\Phi_t,h_t}(x)\|^2_{op}.
\end{equation}
    
By combining Lemma \ref{lemma} and Theorem \ref{wln} we obtain
 \begin{equation}
 \begin{split}
     \mathcal{W}_t(x)&=|\beta_t(x)|^2\sum\limits_{i=1}^{n-1}|a_{i,t}(x)|^{-2} \\
    &=t^2|f(x)|^2\frac{h_{1,t}(x)}{h_{n,t}(x)}\sum\limits_{i=1}^{n-1}\frac{h_{i,t}(x)}{h_{i+1,t}(x)},
    \end{split}
\end{equation}
 which does not depend on the choice of coordinate.
 For $q_n=0$ we have $\beta_t=0$ that implies $C_{\Phi_t,h_t}=0$ from which we have $\mathcal{W}_t=0$.
 \subsubsection{Estimates on asymptotics}
 Next, we present some results from \cite{collier2014asymptotics}  that we need for our proof.
 \begin{theorem}
     For every $x \in \Sigma$ away from the zeros of $q_n$ let $t\rightarrow\infty$, locally in the Hitchin section the metric $h_{h,t}$ on $K^{\frac{n_+1-2j}{2}}$ admits an expansion 
         \begin{equation}\label{h}
             h_{j,t}=(t|q_n|)^{-\frac{n+1-2j}{n}}\big(1+O\big(t^{-\frac{2}{n}}\big)\big)
         \end{equation}
         for all $j$.
 \end{theorem}
 This theorem was originally stated and proved in \cite[Theorem 3.1]{collier2014asymptotics}.

 \begin{proposition}\label{Wt}
     At every $x \in \Sigma$ away from the zeros of $q_n$ the wall function \ref{w} has the following asymptotic behavior:
     \begin{equation}\label{wt}
         \mathcal{W}_t(x)=n-1+O\big(t^{-\frac{2}{n}}\big).
     \end{equation}
     Thus when $t \rightarrow\infty$ we have $\mathcal{W}_t(x)\rightarrow n-1.$
 \end{proposition}
\begin{proof}
  By equations \ref{*} and \ref{h} for each $i$ we obtain
  \begin{equation*}
  \begin{split}
      |a_{i_t}|&=(t|q_n|)^{\frac{1}{n}}\big(1+O\big(t^{-\frac{2}{n}}\big)\big) \\
      |\beta_t|&=(t|q_n|)^{\frac{1}{n}}\big(1+O\big(t^{-\frac{2}{n}}\big)\big).
      \end{split}
  \end{equation*}
  Which implies that $\big|\frac{\beta_t}{a_{i,t}}\big|=1+O\big(t^{-\frac{2}{n}}\big)$.
  By summing the $n-1$ terms of the identity \eqref{id} we obtain the \eqref{wt} which proves the result.
\end{proof}
\subsubsection{The disproof}
Finally, we provide the argument that shows that the main flatness conjecture for $\SL_n(\C)$-Fock bundles cannot be true. More precisely, we will see that there is a positive wall at which universal existence and universal uniqueness cannot both hold in rank $n\geq 3$. In what follows we thus only consider $n\geq3$.
We are now ready to define the positivity wall we use in our proof.

Let 
\begin{equation}
    \Lambda(t)=\max\limits_{x\in\Sigma}\mathcal{W}_t(x)
\end{equation}
which is continuous since $\mathcal{W}_t(x)$ is continuous and $\Sigma$ is compact. When $\mathcal{Q}_0=0$ then $\mathcal{W}_0=0$.
For every $x \in \Sigma$ away from the zeros of $q_n$, Proposition \ref{Wt} implies that $\Lambda(t)\rightarrow n-1>1$ for every sufficiently large $t$ since $n\geq3$. We define 
\begin{equation}
    t_*=\inf\{t>0:\Lambda(t)\geq1\}.
\end{equation}
Then for $0<t_*<\infty$ and $0\leq t< t_*$ we have 
\begin{equation}\label{lambda}
    \Lambda(t)<1, \hspace{0.5cm} \Lambda(t_*)=1.
\end{equation}
We call $t_*$ the first positivity wall. Its main property is that $h_{t_*}$ is not strongly positive while $h_t$ is for $t<t_*$.

\begin{center}
\begin{tikzpicture}[scale=1.7]

    \draw[magenta!90!black, thick]
        plot[smooth cycle, tension=0.8] coordinates {
         (-3.0,0.0)
        (-2.9,0.6)
        (-2.6,1.2)
        (-1.5,1.7)
        (0.2,1.8)
        (2.0,1.3)
        (3.0,0.5)
        (3.0,-0.4)
        (2.3,-0.8)
        (1.2,-0.95)
        (0,-1.0)
        (-1.2,-0.95)
        (-2.0,-0.8)
        (-2.6,-0.5)
        (-2.9,-0.2)
        };

    \draw[pink!90!magenta, line width=4pt, line cap=round]
        plot[smooth, tension=0.8] coordinates {
            (1.8,1.9)
            (1.5,0.9)
            (1.45,0)
            (1.6,-1.5)
        };

    \fill[magenta!90!black] (0.5,0.15) circle (2.5pt);
    \node[above left] at (0.5,0.15) {$t$};
    \node[below left] at (2,0.15) {$h_t$ strongly positive};

    \fill[magenta!90!black] (1.45,0.65) circle (2.5pt);
    \node[above left] at (1.45,0.65) {$t^*$};
    \node[below right] at (1,0.5) {$h_{t_*}$ \textbf{not} strongly positive};

\end{tikzpicture}
\end{center}

\begin{theorem}
   For every rank $n\geq3$ the main flatness conjecture \cite[Conjecture 4.14]{kydonakis2023} is false. Specifically, at least one of the following conditions is forced by the first positivity wall:
   \begin{enumerate}[label=\roman*)]
       \item the field $\Phi_{t_*}$ has no strongly positive flat compatible Hermitian metric $h_{t_*}$,
       \item for all $t<t_*$ sufficiently close to $t_*$, the field $\Phi_t$ has at least two distinct strongly positive flat compatible Hermitian metrics $h_t, \hat{h}_t$.
   \end{enumerate}
\end{theorem}
\begin{proof}
\begin{enumerate}[label=\roman*)]
   
\item  By equation \ref{lambda}, Lemma \ref{lemma} and Proposition \ref{connection}, $h_t$ is strongly positive flat Hermitian metric for every $t<t_*$ but $h_{t_*}$ is not strongly positive. There is a splitting but since it does not satisfy the strongly positivity condition it is not the canonical 3-term connection because Theorem \ref{filling thm} may fail. Since $h_{t_*}$ is no strongly positive flat Hermitian metric, there is no strongly positive metric for $\Phi_{t_*}$ so (a) holds. 
   
\item Assume now there is a strongly positive flat Hermitian metric for $\Phi_{t_*}$ and denote it by $\hat{h}_*$.  By Proposition \ref{local existence}, since $\hat{h}_{*}$ is strongly positive flat for $\Phi_{t_*}$, there is a family of strongly positive flat Hermitian metrics $\hat{h}_t$ for $\Phi_t$ for $t$ close to $t_*$. For $t<t_*$ the metric $h_t$ is another strongly positive flat Hermitian metric for $\Phi_t$. These two families are distinct for every $t<t_*$ sufficiently close to $t_*$. Indeed, let $\{t_j\}_{j\in \N}$ be a sequence such that $t_j<t_*$, and with $t_j\rightarrow t_*$ and assume by way of contradiction that $\hat{h}_{t_j} = h_{t_j}$. Then by smoothness, $\hat{h}_*= \lim\limits_{j\to \infty} \hat{h}_{t_j} =\lim\limits_{j\to\infty}h_{t_j} = h_{t_*}$ which is a contradiction because $\hat{h}_*$ is strongly positive and $h_{t_*}$ is not. Thus, (b) holds.
\end{enumerate}
\end{proof}
Consequently, for every rank $n\geq3$ existence and uniqueness cannot both hold which implies that the main flatness conjecture is false for $\SL_n(\C)$-Fock bundles. The family reaches a finite positivity wall at which either existence of a strongly positive flat Hermitian metric fails or there is a second strongly positive flat Hermitian metric immediately before the wall.
 \begin{remark}
     For rank $n=2$ the norm \eqref{wt} tends to $1$ and thus our approach does not work. This is not an accident since for $n=2$ the main flatness conjecture was proven to hold in \cite{kydonakis2023}.
 \end{remark}

\bibliography{bibliography}{}
\bibliographystyle{amsalpha}

\end{document}